\documentclass[12pt, reqno]{amsart}
\usepackage{graphicx}
\usepackage[hidelinks]{hyperref}
\usepackage{mathrsfs}
\usepackage{amsfonts,amsmath,amsthm,amssymb}
\usepackage{latexsym,euscript,epic,eepic,color}
\usepackage{indentfirst}
\usepackage{tikz}
\usepackage{enumerate}

\usepackage[centering]{geometry}
\usepackage[T1]{fontenc}
\usepackage{lmodern}
\usepackage[expansion=false]{microtype}
\usepackage{amsmath,amsthm,amssymb,mathtools,mathrsfs}
\usepackage{enumerate}
\usepackage{aliascnt}
\usepackage{graphicx,placeins,float}
\usepackage{tikz}
\usepackage[hidelinks]{hyperref}
\usepackage[square,comma,sort&compress,numbers]{natbib}

\theoremstyle{plain}
\newtheorem{theorem}{Theorem}[section]
\newaliascnt{proposition}{theorem}
\newtheorem{proposition}[proposition]{Proposition}
\aliascntresetthe{proposition}
\newaliascnt{lemma}{theorem}
\newtheorem{lemma}[lemma]{Lemma}
\aliascntresetthe{lemma}
\newaliascnt{corollary}{theorem}

\aliascntresetthe{corollary}
\theoremstyle{definition}
\newaliascnt{definition}{theorem}

\aliascntresetthe{definition}
\newaliascnt{example}{theorem}

\aliascntresetthe{example}
\theoremstyle{remark}
\newaliascnt{remark}{theorem}

\aliascntresetthe{remark}

\newcommand{\dimH}{\dim_{\mathrm H}}

\usepackage[nameinlink,capitalize,noabbrev]{cleveref}
\crefname{example}{Example}{Examples}
\Crefname{example}{Example}{Examples}
\numberwithin{equation}{section}

\title[Quantitative multiple recurrence in integer bases]{Hausdorff dimension of
quantitative multiple recurrence sets for $\beta$-transformations I: integer bases and return times $2n$ and $3n$
\texorpdfstring{{}}{exponential targets at times 2n and 3n in integer bases}}
\author{Lingmin Liao$^{1}$}
\address{$^{1}$School of Mathematics and Statistics, Wuhan University, Wuhan, Hubei, China}
\email{lmliao@whu.edu.cn}
\author{Wen Ping$^{1}$}
\address{$^{1}$School of Mathematics and Statistics, Wuhan University, Wuhan, Hubei, China}
\email{pingwen2215@163.com}
\author{Na Yuan$^{2}$}
\address{$^{2}$School of Mathematics, Guangdong University of Education, Xingang Middle Road, 351, Guangzhou, P.~R.\ China}
\email{yuanna199@gmail.com}
\thanks{}
\date{}

\newif\ifshowrevisions
\showrevisionstrue

\begin{document}
\begin{abstract}
Let \(\beta>1\). The $\beta$-transformation on $[0,1)$ is defined by \(T_\beta x=\beta x-\lfloor \beta x\rfloor\). Suppose $\beta=b$ is an integer. 
For fixed \(v_1,v_2\ge0\), 
we  obtain the Hausdorff dimension of the quantitative multiple recurrence set with return times $2n$ and $3n$, i.e., the set of points
$x\in[0,1)$ for which
$ |T_b^{2n}x-x|<b^{-v_1n}$ and
$|T_b^{3n}x-x|<b^{-v_2n}$ hold simultaneously for infinitely
many $n\in\mathbb N$.
In a companion paper, we consider general
non-integer bases $\beta>1$ and study
multiple recurrence sets with return times
$pn$ and $qn$, for arbitrary fixed integers
$1\leq p<q$ and general return rate functions
$\psi_1,\psi_2:\mathbb N\to(0,\infty)$.
\end{abstract}
\maketitle
\thispagestyle{plain}
\pagestyle{plain}

\section{Introduction}\label{sec:introduction}
This paper and its companion \cite{LiaoPingYuanII} concern the Hausdorff dimension of quantitative multiple recurrence sets for $\beta$-transformations. The present paper considers integer bases $\beta=b\geq2$  and the return times \(2n\) and \(3n\), while the companion paper treats arbitrary bases \(\beta>1\) and return times \(pn\) and \(qn\) for each fixed pair of integers \(1\leq p<q\).
\subsection{Background and motivation}
Recurrence is a central topic in dynamical systems. It concerns the
return of an orbit to a neighbourhood of its initial point. 
Let $(X,\rho)$ be a separable metric space, $T\colon X\to X$ a measurable transformation and $\mu$ a finite $T$-invariant Borel measure.
 The Poincar\'e recurrence theorem implies that almost all points are recurrent, i.e.,
\[
\liminf_{n\to\infty}\rho(T^n x,x)=0
\]
for \(\mu\)-almost every \(x\in X\). On the other hand, for a compact topological dynamical
system, Birkhoff's recurrence theorem guarantees the existence of
such recurrent points.

Multiple recurrence concerns the simultaneous return of several iterates of a point to its initial point. 
One form of the problem asks for a point \(x\) such that there exists a sequence $\{n_k\}_{k\geq 1}$ of integers satisfying
\[
 T^{jn_k}x\longrightarrow x,\qquad\text{as}\qquad n_k\rightarrow\infty\qquad\text{for all}\qquad 1\leq j\leq d.
\]
Existence results of such points have been obtained both in measure-theoretic and
topological dynamics; see, for example,
\cite{Furstenberg1977,Furstenberg1981,BergelsonHostKra2005,HuangShaoYe2024}.
However, all these results do not, by themselves, prescribe the recurrence rates. 

The quantitative recurrence studies the recurrent points for which the size of returning neighbourhoods will be prescribed.
Let $(X,\rho)$ be a separable
metric space, $T:X\to X$ a measurable transformation, and
\(\psi:\mathbb N\to(0,\infty)\).
Define the recurrence set by
\[
R(T,\psi)
=
\left\{
x\in X:
\rho(T^n x,x)<\psi(n)
\text{ for infinitely many }n\in\mathbb N
\right\}.
\]
The pioneer work of Boshernitzan \cite{Boshernitzan1993} established the measure of $R(T,\psi)$ for measure-preserving systems. Barreira and Saussol
\cite{Barreira2001} related the recurrence rates of the recurrent points to the local dimension of
the invariant measures.
Further measure and dimension results for recurrence sets
in different dynamical settings can be found in
\cite{Baker2019,Chang2018,He2023QuantitativeRP,Persson2024,Hussain2021,Tan2011,Yuan2023}.

Let $\beta>1$ be a real number. 
We consider the $\beta$-transformation
\[
 T_\beta:[0,1)\longrightarrow[0,1),\qquad
 T_\beta x=\beta x-\lfloor \beta x\rfloor.
\]
For $\psi:\mathbb N\to(0,1]$, let
\[
 R_\beta(\psi)=\{x\in[0,1):|T_\beta^n x-x|<\psi(n)
                   \text{ for infinitely many }n\in\mathbb N\}.
\]
The first quantitive recurrent result of $T_\beta$ was obtained by Tan and Wang \cite{Tan2011} who proved that
\begin{equation}\label{eq-tan-wang}
 \dimH R_\beta(\psi)=\frac1{1+\tau},\quad
 \text{with} \quad\tau=\liminf_{n\to\infty}\frac{-\log_\beta\psi(n)}n,
\end{equation}
where $\dimH$ stands for the Hausdorff dimension.

The purpose of our two companion papers is to study the corresponding multiple
quantitative recurrence problem for \(\beta\)-transformations.
Fix two positive functions
$\psi_1,\psi_2:\mathbb N\to(0,\infty)$. We define the multiple recurrent set
\begin{equation}\label{eq:recurrence-set}
 R_\beta(p,q,\psi_1,\psi_2):=
 \left\{x\in[0,1):
 \begin{array}{l}
 |T_\beta^{pn}x-x|<\psi_1(n),\\
 |T_\beta^{qn}x-x|<\psi_2(n)
 \end{array}
 \ \text{for infinitely many }n\in\mathbb N\right\}.
\end{equation}
In the present paper, we determine the
Hausdorff dimension of
$R_b(2,3,\psi_1,\psi_2)$ for integer bases
$\beta=b\geq2$ and recurrent rate functions
$\psi_i(n)=b^{-v_i n}$, $i=1,2$ with $v_1,v_2\geq0$.
For this case we write $R_b(2,3,v_1,v_2)$ instead of $R_b(2,3,\psi_1,\psi_2)$. 
In the companion paper \cite{LiaoPingYuanII}, we determine
the Hausdorff dimension of
$R_\beta(p,q,\psi_1,\psi_2)$ for real bases $\beta>1$, general positive functions
$\psi_1,\psi_2:\mathbb N\to(0,\infty)$ and arbitrary
integers $1\leq p<q$.

We stress that both inequalities in \eqref{eq:recurrence-set} must hold at the same
integers $n$. To make this requirement explicit, put
\[
E_{1,n}
=\{x\in[0,1):|T_\beta^{pn}x-x|<\psi_1(n)\},\quad
E_{2,n}
=\{x\in[0,1):|T_\beta^{qn}x-x|<\psi_2(n)\}.
\]
Then, 
\begin{equation}\label{eq-same-index-limsup}
 R_\beta(p,q,\psi_1,\psi_2)
 =\limsup_{n\to\infty}(E_{1,n}\cap E_{2,n}),
\end{equation}
whereas
\begin{equation}\label{eq-limsup-intersection-inclusion}
 \limsup_{n\to\infty}(E_{1,n}\cap E_{2,n})
 \subseteq
 \left(\limsup_{n\to\infty}E_{1,n}\right)
 \cap\left(\limsup_{n\to\infty}E_{2,n}\right).
\end{equation}

Remark that the two infinite sets of return integers on the right may be disjoint;
for general sequences of sets, the reverse inclusion is false.
Accordingly, a dimension estimate for the right-hand side alone does
not determine the dimension on the left.
This also distinguishes the present problem from an application of
closure under intersections in Falconer's large-intersection classes
\cite{Falconer1994}. Note that large-intersection results for recurrence sets of
hyperbolic toral automorphisms were obtained in \cite{Persson2024}. However, intersecting two such limsup sets
does not impose common return integers.

On the other hand, a direct use of a mass
transference theorem, such as \cite[Theorem 1.4]{He2025MTP}, would also
require the measure and size hypotheses to be verified for sets inside
the simultaneous events $E_{1,n}\cap E_{2,n}$, and will be delicate. 
Instead, we estimate the
simultaneous events $E_{1,n}$ and $ E_{2,n}$ directly.
\subsection{Main results and examples}
\label{subsec:main-results}
\leavevmode\par
Now, we state our main result of the present paper.
\begin{theorem}\label{thm:exp}
For every integer $b\geq2$ and every $v_1,v_2\geq0$, we have
\begin{equation}\label{eq:exp-D}
 \dimH R_b(2,3,v_1,v_2)=D_{2,3}(v_1,v_2):=
 \begin{cases}
 \displaystyle\min\left\{\frac{2}{2+v_1},\ \frac{3-v_1}{3+v_2}\right\},
       &0\leq v_1\leq1,\\[11pt]
 \displaystyle\frac{\max\{1,\ 3-v_1,\ 2-v_2\}}{\max\{2+v_1,\ 3+v_2\}},
       &v_1>1.
 \end{cases}
\end{equation}
\end{theorem}

Let us give some ideas on how such dimension formula comes about. 
First to illustrate how the threshold $v_1=1$ arises,
consider the case $b=2$.
Fix an integer $n\geq1$ and assume that
$v_1n$ and $v_2n$ are positive integers.
For $x\in[0,1)$, write
\[
x=\sum_{i=1}^{\infty}\varepsilon_i2^{-i},
\qquad \varepsilon_i\in\{0,1\},
\]
using the expansion ending in zeros whenever
$x$ has two binary expansions.
The binary digits of $T_2^{2n}x$ are
$\varepsilon_{2n+1},\varepsilon_{2n+2},\ldots$.
Hence, the equalities
\begin{equation}\label{2n}
\varepsilon_{2n+i}=\varepsilon_i
\qquad(1\leq i\leq v_1n)
\end{equation}
imply both $x$ and $T_2^{2n}x$ in the interval
\([
\sum_{i=1}^{v_1n}\varepsilon_i2^{-i},
\ \sum_{i=1}^{v_1n}\varepsilon_i2^{-i}
   +2^{-v_1n}
).
\)
Consequently,
\(|T_2^{2n}x-x|<2^{-v_1n}.
\)
Similarly, the equalities
\begin{equation}\label{3n}
\varepsilon_{3n+i}=\varepsilon_i
\qquad(1\leq i\leq v_2n)
\end{equation}
imply
\(|T_2^{3n}x-x|<2^{-v_2n}.
\)
Thus \eqref{2n} and \eqref{3n} are sufficient
for the two return inequalities defining
$R_2(2,3,v_1,v_2)$ to hold at the fixed
integer $n$.

However, equation~\eqref{2n} prescribes the digits
at positions $2n+1,\ldots,2n+v_1n$,
while \eqref{3n} prescribes the digits
at positions $3n+1,\ldots,3n+v_2n$.
When $v_1\leq1$, the prescribed position
ranges are disjoint, since
$2n+v_1n\leq3n$.
When $v_1>1$, both \eqref{2n} and \eqref{3n}
prescribe digits at the positions
\(3n+1,\ldots,\min\{2n+v_1n,3n+v_2n\}.
\)
In fact, for every integer $j$ satisfying
\(1\leq j\leq\min\{(v_1-1)n,v_2n\},
\)
we have $n+j\leq v_1n$ and $j\leq v_2n$.
Applying \eqref{2n} with $i=n+j$ and
\eqref{3n} with $i=j$ gives
\[
\varepsilon_{n+j}
=\varepsilon_{3n+j}
=\varepsilon_j.
\]
Therefore, 
\(v_1=1\) is the threshold beyond which the two prescribed
position ranges of \eqref{2n}
and \eqref{3n} begin to overlap.

The effect of the imposed digit equalities can be
illustrated by the following two choices of exponents.

We first compare the cases $v_1=1$ and
$v_1=2$, keeping $v_2=1$.

First, let $v_1=v_2=1$.
Equations~\eqref{2n} and~\eqref{3n}
determine the digits at positions
$2n+1,\ldots,3n$ and $3n+1,\ldots,4n$,
respectively, by copying the first $n$ digits.
Therefore,
the first $2n$ digits can  be
chosen freely and uniquely determine
the first $4n$ digits.
The proportion of freely prescribable digits
is $2n/(4n)=1/2$, and \eqref{eq:exp-D} gives
\[
D_{2,3}(1,1)=\frac12,
\]
which coincide our intuition. 

Now let $v_1=2$ and $v_2=1$.
For $1\leq j\leq n$, applying \eqref{2n}
with $i=n+j$ and \eqref{3n} with $i=j$ gives
\[
\varepsilon_{n+j}
=\varepsilon_{3n+j}
=\varepsilon_j.
\]
Thus, the first $n$ digits determine
the digits at positions $n+1,\ldots,2n$
and $3n+1,\ldots,4n$.
Further, the imposed equalities \(\varepsilon_{2n+i}=\varepsilon_i\)  \((1\leq i\leq 2n)\) 
determine all digits at positions $2n+1,\ldots,4n$.
Thus, the first $4n$ digits are completely
determined by the first $n$ digits, and
\(D_{2,3}(2,1)=\frac14.
\)

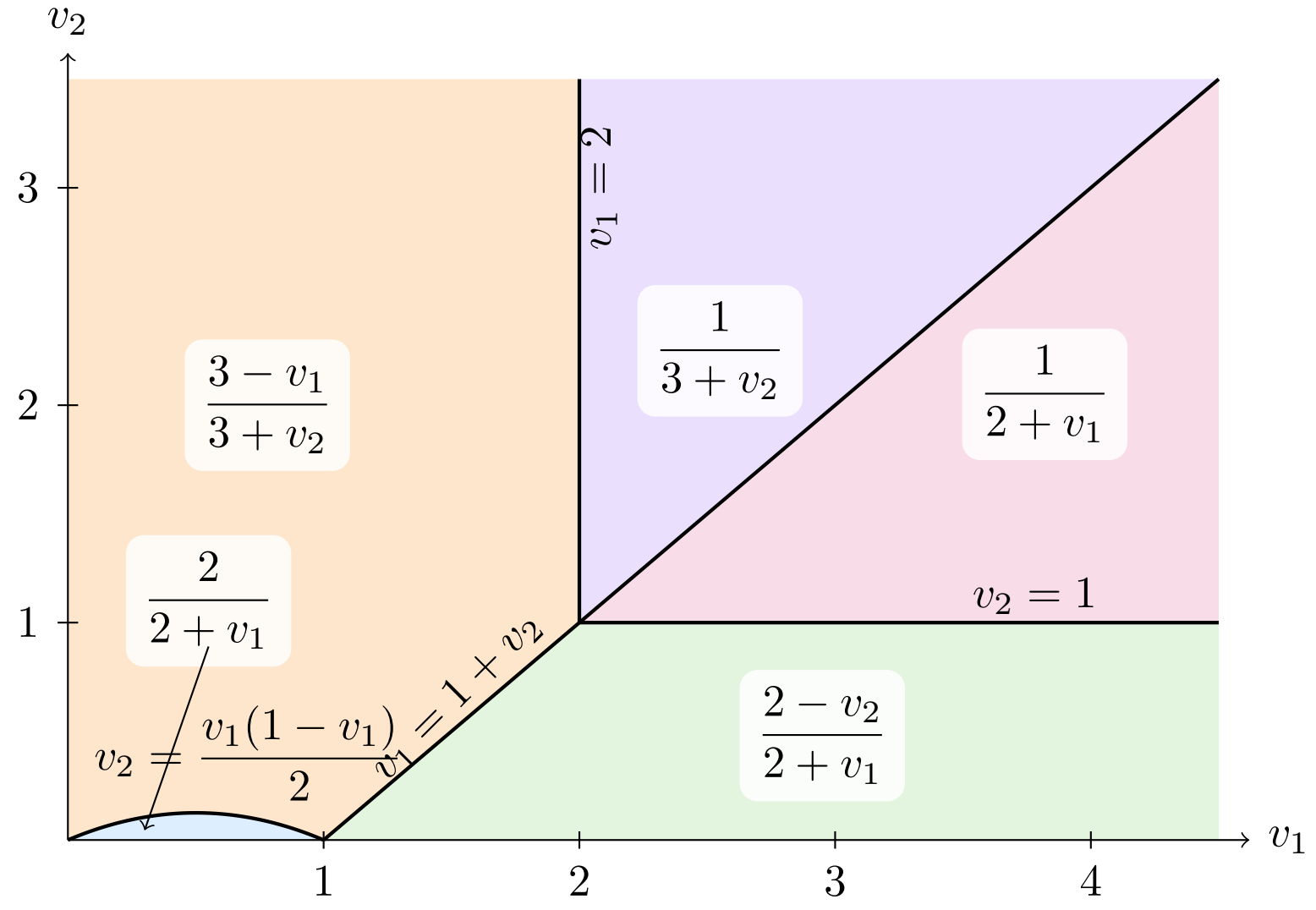
\begin{figure}[htbp]
 \centering
 \begin{tikzpicture}[x=2.0cm,y=1.7cm,font=\footnotesize]
  \definecolor{regblue}{RGB}{220,238,255}
  \definecolor{regorange}{RGB}{253,230,204}
  \definecolor{reggreen}{RGB}{227,244,223}
  \definecolor{regpurple}{RGB}{234,223,253}
  \definecolor{regrose}{RGB}{248,220,231}
  \clip (-0.48,-0.40) rectangle (4.95,3.98);
  \fill[regorange]
    (0,0) plot[domain=0:1,samples=80] (\x,{\x*(1-\x)/2})
    -- (2,1) -- (2,3.5) -- (0,3.5) -- cycle;
  \fill[regblue]
    (0,0) plot[domain=0:1,samples=80] (\x,{\x*(1-\x)/2})
    -- cycle;
  \fill[reggreen] (1,0) -- (2,1) -- (4.5,1) -- (4.5,0) -- cycle;
  \fill[regpurple] (2,1) -- (2,3.5) -- (4.5,3.5) -- cycle;
  \fill[regrose] (2,1) -- (4.5,3.5) -- (4.5,1) -- cycle;
  \draw[thick]
    plot[domain=0:1,samples=80] (\x,{\x*(1-\x)/2});
  \draw[thick] (1,0) -- (2,1) -- (4.5,1);
  \draw[thick] (2,1) -- (4.5,3.5);
  \draw[thick] (2,1) -- (2,3.5);
  \draw[->] (0,0) -- (4.62,0) node[right] {$v_1$};
  \draw[->] (0,0) -- (0,3.62) node[above] {$v_2$};
  \foreach \x in {1,2,3,4}
    \draw (\x,0.04) -- (\x,-0.04) node[below] {\footnotesize\(\x\)};
  \foreach \y in {1,2,3}
    \draw (0.04,\y) -- (-0.04,\y) node[left] {\footnotesize\(\y\)};
  \node[fill=white,fill opacity=.86,text opacity=1,rounded corners]
    at (.55,1.1) {$\dfrac{2}{2+v_1}$};
  \draw[->] (.55,.89) -- (.30,.045);
  \node[fill=white,fill opacity=.86,text opacity=1,rounded corners]
    at (.78,2) {$\dfrac{3-v_1}{3+v_2}$};
  \node[fill=white,fill opacity=.86,text opacity=1,rounded corners]
    at (2.95,.48) {$\dfrac{2-v_2}{2+v_1}$};
  \node[fill=white,fill opacity=.86,text opacity=1,rounded corners]
    at (2.55,2.25) {$\dfrac{1}{3+v_2}$};
  \node[fill=white,fill opacity=.86,text opacity=1,rounded corners]
    at (3.82,2.05) {$\dfrac{1}{2+v_1}$};
  \node[anchor=west,font=\footnotesize] at (.03,0.4)
    {$v_2=\dfrac{v_1(1-v_1)}2$};
  \node[rotate=45] at (1.52,.66) {$v_1=1+v_2$};
  \node at (3.78,1.12) {$v_2=1$};
  \node[rotate=90] at (2.08,3.0) {$v_1=2$};
 \end{tikzpicture}
\caption{Partition of the $(v_1,v_2)$-plane according to the
formula for $D_{2,3}(v_1,v_2)$. The expression written in each
colored region is the value of $D_{2,3}$ for parameters in that
region. The values are equal on the boundary between any two regions, and the mapping $(v_1,v_2)\mapsto D_{2,3}(v_1,v_2)$ is continuous. Only a bounded part of the
parameter plane is displayed.}
 \label{fig:23-regions}
\end{figure}

Figure~\ref{fig:23-regions} displays a partition
of the $(v_1,v_2)$-plane according to the formula
for $D_{2,3}(v_1,v_2).$
By \eqref{eq:exp-D}, 
\[
 D_{2,3}(v_1,v_2)=
 \begin{cases}
\dfrac{2}{2+v_1},
& 0\le v_1\le 1,\quad
0\le v_2\le \frac{v_1(1-v_1)}{2},
\\[6pt]
 \dfrac{3-v_1}{3+v_2},
   &0 \le v_1 \le 2,\quad
v_2 \ge \max\left\{v_1-1,\frac{v_1(1-v_1)}{2}\right\},\\[6pt]
 \dfrac{2-v_2}{2+v_1},
   &0\leq v_2\leq1,\quad v_1\geq1+v_2,\\[6pt]
 \dfrac1{3+v_2},
   &v_2\geq1,\quad 2\leq v_1\leq1+v_2,\\[6pt]
 \dfrac1{2+v_1},
   &v_2\geq1,\quad v_1\geq1+v_2.\\[6pt]
 \end{cases}
\]

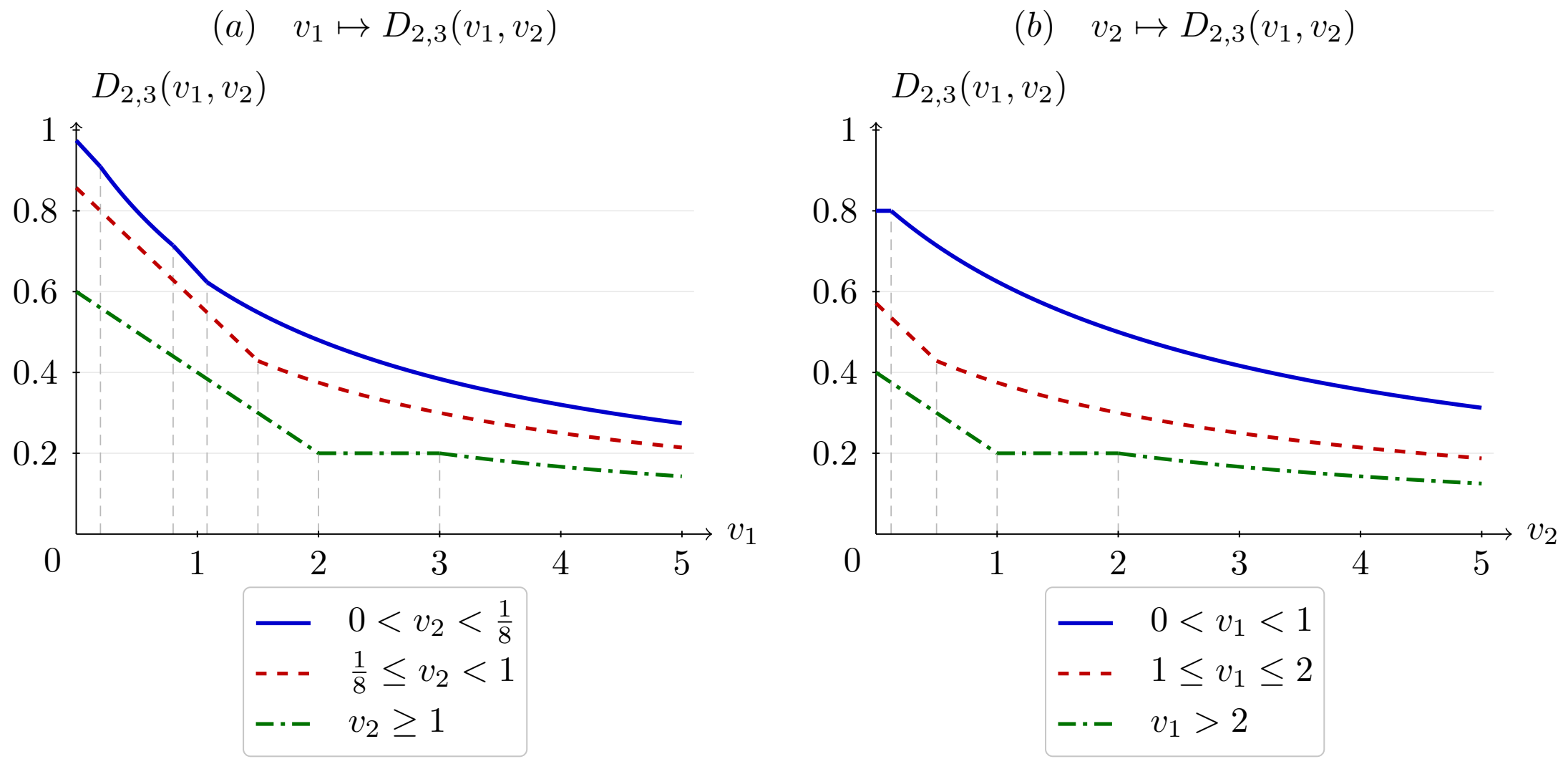
\begin{figure}[htbp]
\centering

\tikzset{
  dA/.style={blue!80!black,line width=1.15pt},
  dB/.style={red!75!black,line width=1.05pt,dashed},
  dC/.style={green!45!black,line width=1.05pt,
    dash pattern=on 5pt off 2pt on 1.2pt off 2pt},
  dGuide/.style={gray!55,densely dashed,line width=.35pt},
  dGrid/.style={gray!18,line width=.3pt}
}

\begin{minipage}[t]{0.47\linewidth}
\centering
\vspace{0pt}
\resizebox{\linewidth}{!}{%
\begin{tikzpicture}[x=1.20cm,y=4.0cm,font=\footnotesize]

  \node[anchor=south] at (2.55,1.17)
    {$(a)\quad v_1\mapsto D_{2,3}(v_1,v_2)$};

  \foreach \y in {0.2,0.4,0.6,0.8}
    \draw[dGrid] (0,\y)--(5.1,\y);

  \draw[dGuide] (0.2,0)--(0.2,{2/2.2});
  \draw[dGuide] (0.8,0)--(0.8,{2/2.8});
  \draw[dGuide] (1.08,0)--(1.08,{1.92/3.08});
  \draw[dGuide] (1.5,0)--(1.5,{1.5/3.5});
  \draw[dGuide] (2,0)--(2,0.2);
  \draw[dGuide] (3,0)--(3,0.2);

  \draw[->] (0,0)--(5.25,0)
    node[right] {$v_1$};
  \draw[->] (0,0)--(0,1.02)
    node[above right] {$D_{2,3}(v_1,v_2)$};

  \foreach \x in {1,2,3,4,5}
    \draw (\x,0.008)--(\x,-0.008)
      node[below] {$\x$};

  \foreach \y/\lab in
    {0.2/0.2,0.4/0.4,0.6/0.6,0.8/0.8,1/1}
    \draw (0.03,\y)--(-0.03,\y)
      node[left] {$\lab$};

  \node[below left] at (0,0) {$0$};

  \draw[dA]
    plot[domain=0:0.2,samples=2]
      (\x,{(3-\x)/3.08});
  \draw[dA]
    plot[domain=0.2:0.8,samples=90]
      (\x,{2/(2+\x)});
  \draw[dA]
    plot[domain=0.8:1.08,samples=2]
      (\x,{(3-\x)/3.08});
  \draw[dA]
    plot[domain=1.08:5,samples=120]
      (\x,{1.92/(2+\x)});

  \draw[dB]
    plot[domain=0:1.5,samples=2]
      (\x,{(3-\x)/3.5});
  \draw[dB]
    plot[domain=1.5:5,samples=120]
      (\x,{1.5/(2+\x)});

  \draw[dC]
    plot[domain=0:2,samples=2]
      (\x,{(3-\x)/5});
  \draw[dC] (2,0.2)--(3,0.2);
  \draw[dC]
    plot[domain=3:5,samples=120]
      (\x,{1/(2+\x)});

  \node[
    anchor=north,
    fill=white,
    rounded corners=2pt,
    draw=gray!45,
    inner sep=3pt
  ] at (2.55,-0.13)
  {\begin{tabular}{@{}ll@{}}
    \tikz[baseline=-0.5ex]{
      \draw[dA] (0,0)--(0.45,0);
    }
    & $\vphantom{\frac18}0<v_2<\frac18$ \\[1pt]
    \tikz[baseline=-0.5ex]{
      \draw[dB] (0,0)--(0.45,0);
    }
    & $\vphantom{\frac18}\frac18\leq v_2<1$ \\[1pt]
    \tikz[baseline=-0.5ex]{
      \draw[dC] (0,0)--(0.45,0);
    }
    & $\vphantom{\frac18}v_2\geq1$
  \end{tabular}};

\end{tikzpicture}%
}
\end{minipage}%
\hspace{0.015\linewidth}%
\begin{minipage}[t]{0.47\linewidth}
\centering
\vspace{0pt}
\resizebox{\linewidth}{!}{%
\begin{tikzpicture}[x=1.20cm,y=4.0cm,font=\footnotesize]

  \node[anchor=south] at (2.55,1.17)
    {$(b)\quad v_2\mapsto D_{2,3}(v_1,v_2)$};

  \foreach \y in {0.2,0.4,0.6,0.8}
    \draw[dGrid] (0,\y)--(5.1,\y);

  \draw[dGuide] (0.125,0)--(0.125,0.8);
  \draw[dGuide] (0.5,0)--(0.5,{1.5/3.5});
  \draw[dGuide] (1,0)--(1,0.2);
  \draw[dGuide] (2,0)--(2,0.2);

  \draw[->] (0,0)--(5.25,0)
    node[right] {$v_2$};
  \draw[->] (0,0)--(0,1.02)
    node[above right] {$D_{2,3}(v_1,v_2)$};

  \foreach \x in {1,2,3,4,5}
    \draw (\x,0.008)--(\x,-0.008)
      node[below] {$\x$};

  \foreach \y/\lab in
    {0.2/0.2,0.4/0.4,0.6/0.6,0.8/0.8,1/1}
    \draw (0.03,\y)--(-0.03,\y)
      node[left] {$\lab$};

  \node[below left] at (0,0) {$0$};

  \draw[dA] (0,0.8)--(0.125,0.8);
  \draw[dA]
    plot[domain=0.125:5,samples=120]
      (\x,{2.5/(3+\x)});

  \draw[dB]
    plot[domain=0:0.5,samples=2]
      (\x,{(2-\x)/3.5});
  \draw[dB]
    plot[domain=0.5:5,samples=120]
      (\x,{1.5/(3+\x)});

  \draw[dC]
    plot[domain=0:1,samples=2]
      (\x,{(2-\x)/5});
  \draw[dC] (1,0.2)--(2,0.2);
  \draw[dC]
    plot[domain=2:5,samples=120]
      (\x,{1/(3+\x)});

  \node[
    anchor=north,
    fill=white,
    rounded corners=2pt,
    draw=gray!45,
    inner sep=3pt
  ] at (2.55,-0.13)
  {\begin{tabular}{@{}ll@{}}
    \tikz[baseline=-0.5ex]{
      \draw[dA] (0,0)--(0.45,0);
    }
    & $\vphantom{\frac18}0<v_1<1$ \\[1pt]
    \tikz[baseline=-0.5ex]{
      \draw[dB] (0,0)--(0.45,0);
    }
    & $\vphantom{\frac18}1\leq v_1\leq2$ \\[1pt]
    \tikz[baseline=-0.5ex]{
      \draw[dC] (0,0)--(0.45,0);
    }
    & $\vphantom{\frac18}v_1>2$
  \end{tabular}};

\end{tikzpicture}%
}
\end{minipage}
\par\smallskip

\caption{Graphs of $D_{2,3}(v_1,v_2)$ with
(a) $v_2$ fixed and (b) $v_1$ fixed.
The solid blue, dashed red, and dash-dotted green curves
correspond, respectively, to $v_2=0.08,0.5,2$ in (a)
and $v_1=0.5,1.5,3$ in (b).}
\label{fig:23-sections-combined}
\end{figure}

The formulas agree wherever two of the displayed regions share a
boundary.
Figure~\ref{fig:23-sections-combined}(a) shows the graphs of
$v_1\mapsto D_{2,3}(v_1,v_2)$ for the fixed values
$v_2=0.08$, $0.5$, and $2$, which belong to the ranges
$0<v_2<1/8$, $1/8\leq v_2<1$, and $v_2\geq1$, respectively.
Figure~\ref{fig:23-sections-combined}(b) shows the graphs of
$v_2\mapsto D_{2,3}(v_1,v_2)$ for the fixed values
$v_1=0.5$, $1.5$, and $3$, which belong to the ranges
$0<v_1<1$, $1\leq v_1\leq2$, and $v_1>2$, respectively.

\subsection{Organization}
Section~\ref{sec:preliminaries} recalls the basic properties of \(b\)-ary expansions used in the proofs. 
Section~\ref{sec:exp-upper} gives the upper bound of Theorem~\ref{thm:exp}. Section~\ref{sec:exp-digits} proves the lower bound.
Our companion paper \cite{LiaoPingYuanII} treats arbitrary bases \(\beta>1\) with return times \(pn\) and \(qn\)\((1\leq p<q)\), and more general target functions. The proof in \cite{LiaoPingYuanII} is different to that in this paper.

\section{Preliminaries}\label{sec:preliminaries}

This section collects some basic facts that
will be used below. Throughout the paper, $\#E$ denotes the cardinality of a finite
set $E$.

Let $b\geq2$ be an integer. The
$b$-transformation $T_b:[0,1)\to[0,1)$ is
defined by
\[
T_bx=bx-\lfloor bx\rfloor,
\]
where $\lfloor\xi\rfloor$ denotes the greatest
integer not exceeding $\xi$.

For $x\in[0,1)$, define
\begin{equation}\label{def}
\varepsilon_n(x,b)
=
\lfloor bT_b^{n-1}x\rfloor,
\qquad n\geq1.
\end{equation}
Then, $x$ has the following  $b$-ary expansion:
\begin{equation}\label{eq-beta-expansion}
x=\sum_{n=1}^{\infty}
\frac{\varepsilon_n(x,b)}{b^n}.
\end{equation}
We also call the sequence
\[
\varepsilon(x,b)
=
(\varepsilon_1(x,b),\varepsilon_2(x,b),\ldots)
\]
the $b$-expansion of $x$.
For every integer $\ell\geq1$, we have
\begin{equation}\label{eq:expansion-finite}
x=
\sum_{j=1}^{\ell}
\frac{\varepsilon_j(x,b)}{b^j}
+\frac{T_b^\ell x}{b^\ell},
\end{equation}
and for any $j\geq1$, $\varepsilon_j(T_b^\ell x,b)
=
\varepsilon_{\ell+j}(x,b)$.
Each digit $\varepsilon_n(x,b)$ belongs to
\begin{equation}\label{eq-A-b}
\mathcal A_b=\{0,1,\ldots,b-1\}.   
\end{equation}

For a word
$\mathbf w=(w_1,\ldots,w_\ell)\in\mathcal A_b^\ell$,
define the order-$\ell$ cylinder associated with
$\mathbf w$ by
\begin{equation}\label{eq-cylinder-def}
I_\ell(\mathbf w)
=
\left\{
x\in[0,1):
\varepsilon_j(x,b)=w_j,\ 
1\leq j\leq\ell
\right\}.
\end{equation}
For $x\in[0,1)$ and $\ell\geq1$, we write
\[
I_\ell(x)
=
I_\ell\bigl(
\varepsilon_1(x,b),\ldots,\varepsilon_\ell(x,b)
\bigr)
\]
for the order-$\ell$ cylinder containing $x$.
The cylinder $I_\ell(\mathbf w)$ is a
half-open interval: 
\[
I_\ell(\mathbf w)
=
\left[
\sum_{j=1}^{\ell}w_jb^{-j},
\sum_{j=1}^{\ell}w_jb^{-j}+b^{-\ell}
\right)
\]
of length $b^{-\ell}$.
At each fixed order $\ell$, the cylinders are
pairwise disjoint and partition $[0,1)$.

For $s>0$, $\rho>0$ and $E\subset\mathbb R$,
define
\[
\mathcal H_\rho^s(E)
=
\inf\left\{
\sum_{i=1}^{\infty}
(\operatorname{diam}U_i)^s:
E\subseteq\bigcup_{i=1}^{\infty}U_i,\ 
\operatorname{diam}U_i\leq\rho
\right\}.
\]
The $s$-dimensional Hausdorff measure of $E$
is defined by
\[
\mathcal H^s(E)
=
\lim_{\rho\to0}\mathcal H_\rho^s(E).
\]
The Hausdorff dimension of a set \(E\) is defined by
\[
\dim_{\mathrm H} E
= \inf \bigl\{ s>0 : \mathcal H^s(E)=0 \bigr\}.
\]
In particular, \(\mathcal H^s(E)=0\) if \(s>\dim_{\mathrm H}E\),
and \(\mathcal H^s(E)=\infty\) if \(0<s<\dim_{\mathrm H}E\).

The following proposition, called the mass distribution principle, is a useful tool to detemine the Hausdorff dimension of a set.
\begin{proposition}[\cite{Falconer1990}]
\label{mmd}
Let $F\subset\mathbb R$ be a Borel set, and
$\mu$ be a Borel probability measure satisfying
$\mu(F)=1$.
Suppose that there exist constants $C>0$,
$s>0$ and $r_0>0$ such that
\[
\mu(J)\leq C|J|^s
\]
for every interval $J\subset\mathbb R$
with $|J|\leq r_0$.
Then
\(\dim_{\mathrm H}F\geq s.
\)
\end{proposition}

We refer to \cite{Falconer1990} for further
background on Hausdorff measure and dimension.

\section{The upper bound}
\label{sec:exp-upper}
The proof of Theorem~\ref{thm:exp} is naturally divides into two parts: the upper bound and the lower bound. In this section, we prove the upper bound. 

Recall 
\[
E_{1,n}
=\{x\in[0,1):|T_b^{2n}x-x|<b^{-v_1n}\}
\quad\text{
and}
\quad
E_{2,n}
=\{x\in[0,1):|T_b^{3n}x-x|<b^{-v_2n}\}.
\]
Put
\[
E_n=E_{1,n}\cap E_{2,n}.
\]
Then
\[
R_b(2,3,v_1,v_2)=\limsup_{n\to\infty}E_n.
\]

For $\mathbf w\in\mathcal{A}_b^{2n}$ and
$\mathbf z\in\mathcal{A}_b^{3n}$, define
\[
J_{2n,1}(\mathbf w)
=I_{2n}(\mathbf w)\cap E_{1,n}\quad\text{and}
\quad
J_{3n,2}(\mathbf z)
=I_{3n}(\mathbf z)\cap E_{2,n}.
\]
For $\mathbf w=(w_1,\ldots,w_{2n})
$ and
$\mathbf z=(z_1,\ldots,z_{3n})$, put
\[
a=\sum_{i=1}^{2n}w_i b^{2n-i},
\qquad
c=\sum_{i=1}^{3n}z_i b^{3n-i}.
\]
By \eqref{eq:expansion-finite},
for any $x\in I_{2n}(\mathbf w)$, we have
\[
|T_b^{2n}x-x|
=
\left|(b^{2n}-1)x-a\right|.
\]
Since $b^{2n}-1>0$, the definition of
$J_{2n,1}(\mathbf w)$ gives
\begin{equation}\label{eq:exp-J2n-condition}
x\in J_{2n,1}(\mathbf w)
\quad\Longleftrightarrow\quad
\frac{a-b^{-v_1n}}{b^{2n}-1}
<x<
\frac{a+b^{-v_1n}}{b^{2n}-1}.
\end{equation}
Since $I_{2n}(\mathbf w)$ is an interval,
$J_{2n,1}(\mathbf w)$ is also an interval.
The bounds for $x$ in
\eqref{eq:exp-J2n-condition} imply
\begin{equation}\label{eq:exp-J2n}
|J_{2n,1}(\mathbf w)|
\leq
\frac{2b^{-v_1n}}{b^{2n}-1}
\leq
4b^{-(2+v_1)n},
\end{equation}
where the last inequalities follows from
$ b^{2n}-1\geq b^{2n}/2$.
Similarly, for $x\in I_{3n}(\mathbf z)$,
the identity
$T_b^{3n}x-x=(b^{3n}-1)x-c$ gives
\begin{equation}\label{eq:exp-J3n-condition}
x\in J_{3n,2}(\mathbf z)
\quad\Longleftrightarrow\quad
\frac{c-b^{-v_2n}}{b^{3n}-1}
<x<
\frac{c+b^{-v_2n}}{b^{3n}-1}.
\end{equation}
Hence $J_{3n,2}(\mathbf z)$ is an interval
and
\begin{equation}\label{eq:exp-J3n}
|J_{3n,2}(\mathbf z)|
\leq
\frac{2b^{-v_2n}}{b^{3n}-1}
\leq
4b^{-(3+v_2)n},
\end{equation}
where the last inequalities 
follows from
$b^{3n}-1\geq b^{3n}/2$.

We first use only the condition
$|T_b^{2n}x-x|<b^{-v_1n}$.
The family
\[
\{J_{2n,1}(\mathbf w):
\mathbf w\in\mathcal A_b^{2n}\}
\]
covers $E_{1,n}$ and hence $E_n$.
Since $\#\mathcal A_b^{2n}=b^{2n}$, 
the formula 
\eqref{eq:exp-J2n} gives, for every $s>0$,
\begin{equation}\label{eq:exp-first-cover}
\sum_{\mathbf w\in\mathcal A_b^{2n}}
|J_{2n,1}(\mathbf w)|^s
\leq
4^s b^{[2-s(2+v_1)]n}.
\end{equation}

We next use both return conditions.
For $\mathbf w\in\mathcal A_b^{2n}$, put
\[
\mathcal V_n(\mathbf w)
:=
\left\{
\mathbf z\in\mathcal A_b^{3n}:
J_{2n,1}(\mathbf w)
\cap J_{3n,2}(\mathbf z)\ne\varnothing
\right\}.
\]
Then, the family
\begin{equation}\label{eq:exp-second-cover}
\left\{
J_{2n,1}(\mathbf w)
\cap J_{3n,2}(\mathbf z):
\mathbf w\in\mathcal A_b^{2n},\
\mathbf z\in\mathcal V_n(\mathbf w)
\right\}
\end{equation}
covers $E_n$.

We now estimate the number of nonempty intersections.
The maps $\mathbf w\mapsto a$ from $\mathcal A_b^{2n}$ to $\{0,\ldots,b^{2n}-1\}$ and
$\mathbf z\mapsto c$ from $\mathcal A_b^{3n}$ to $\{0,\ldots,b^{3n}-1\}$ are bijections
respectively.

Suppose that $\mathbf z\in\mathcal V_n(\mathbf w)$,
and choose
\[
x\in J_{2n,1}(\mathbf w)
\cap J_{3n,2}(\mathbf z).
\]
By \eqref{eq:exp-J2n-condition}, \eqref{eq:exp-J3n-condition}
and the triangle inequality, we have
\[
\begin{aligned}
\left|
\frac{a}{b^{2n}-1}
-\frac{c}{b^{3n}-1}
\right|
&\leq
\left|x-\frac{a}{b^{2n}-1}\right|
+
\left|x-\frac{c}{b^{3n}-1}\right|\\
&<
\frac{b^{-v_1n}}{b^{2n}-1}
+
\frac{b^{-v_2n}}{b^{3n}-1}.
\end{aligned}
\]
Multiplying by
\(\frac{(b^{2n}-1)(b^{3n}-1)}{b^n-1}>0
\)
and using
\[
b^{2n}-1=(b^n-1)(b^n+1),
\qquad
b^{3n}-1=(b^n-1)(b^{2n}+b^n+1),
\]
we obtain
\begin{equation}\label{eq:exp-integer-difference}
\begin{aligned}
\left|
a(b^{2n}+b^n+1)-c(b^n+1)
\right|
<
(b^{2n}+b^n+1)b^{-v_1n}
+(b^n+1)b^{-v_2n}.
\end{aligned}
\end{equation}
The left-hand side of \eqref{eq:exp-integer-difference}
is the absolute value of an integer denoted by $t$. 
Moreover,
the right-hand side is bounded by
\(C_0(
b^{(2-v_1)n}+b^{(1-v_2)n}
),
\)
where $C_0$ is independent of $n$.
Since an interval of radius $R$ contains at most
$2R+1$ integers, the integer $t$ of the left-hand side of
\eqref{eq:exp-integer-difference} can take at most
\[
C\left(
1+b^{(2-v_1)n}+b^{(1-v_2)n}
\right)
\]
different values.

Since $c$ is uniquely determined by $t$ and $a$, once $t$ is determined, we need only to count the possibility of $a$. 
We next count for a fixed integer value $t$, the possible $a$ (on the pairs $(a,c)$) such that $a(b^{2n}+b^n+1)-c(b^n+1)=t$. 
By
\[
b^{2n}+b^n+1=b^n(b^n+1)+1,
\]
we have
\[
a\equiv
a(b^{2n}+b^n+1)-c(b^n+1)= t
\pmod{b^n+1}.
\]
Thus, the residue of $a$ modulo $b^n+1$
is fixed. As
\[
0\leq a\leq b^{2n}-1
=(b^n-1)(b^n+1),
\]
there are at most $b^n$ possible values 
of $a$ for each fixed valued 
$t$. 

	Every nonempty intersection
	$J_{2n,1}(\mathbf w)\cap J_{3n,2}(\mathbf z)$
	gives an integer pair $(a,c)$ satisfying
	\eqref{eq:exp-integer-difference}.
	Since the maps from $\mathbf w$ to $a$
	and from $\mathbf z$ to $c$ are bijections,
	the number of the nonempty intersections in
\eqref{eq:exp-second-cover}
	is bounded by the number of
	integer pairs $(a,c)$ satisfying
	\eqref{eq:exp-integer-difference}.
It follows that we have at most
\begin{equation}\label{eq:exp-intersection-count}
C\left(
1+b^{(2-v_1)n}+b^{(1-v_2)n}
\right)b^n=C\bigl(
b^n+b^{(3-v_1)n}+b^{(2-v_2)n}
\bigr)
\leq
3C b^{n\max\{1,3-v_1,2-v_2\}}
\end{equation}
nonempty intersections $J_{2n,1}(\mathbf w)\cap J_{3n,2}(\mathbf z)$. 
Moreover, every nonempty intersection in
\eqref{eq:exp-second-cover} is contained in
both $J_{2n,1}(\mathbf w)$ and
$J_{3n,2}(\mathbf z)$. By
\eqref{eq:exp-J2n} and \eqref{eq:exp-J3n},
the diameter of $J_{2n,1}(\mathbf w)
\cap J_{3n,2}(\mathbf z)$ is at most
\[
4b^{-n\max\{2+v_1,3+v_2\}}.
\]
Combining this diameter estimate with
\eqref{eq:exp-intersection-count}, for $s>0$, we have
\begin{equation}\label{eq:exp-second-cover-cost}
\sum_{\mathbf w\in\mathcal A_b^{2n}}
\sum_{\mathbf z\in\mathcal V_n(\mathbf w)}
\left|
J_{2n,1}(\mathbf w)
\cap J_{3n,2}(\mathbf z)
\right|^s
\leq
3C\,4^s
b^{n[
\max\{1,3-v_1,2-v_2\}
-s\max\{2+v_1,3+v_2\}
]}.
\end{equation}

Since
$R_b(2,3,v_1,v_2)=\limsup_{n\to\infty}E_n$,
for every integer $N\geq1$, we have
\[
R_b(2,3,v_1,v_2)
\subseteq
\bigcup_{n\geq N}
\ \bigcup_{\mathbf w\in\mathcal A_b^{2n}}
J_{2n,1}(\mathbf w)
\]
and
\[
R_b(2,3,v_1,v_2)
\subseteq
\bigcup_{n\geq N}
\ \bigcup_{\mathbf w\in\mathcal A_b^{2n}}
\ \bigcup_{\mathbf z\in\mathcal V_n(\mathbf w)}
\bigl(
J_{2n,1}(\mathbf w)\cap
J_{3n,2}(\mathbf z)
\bigr).
\]
For $n\geq N$, \eqref{eq:exp-J2n} gives
\[
\bigl|
J_{2n,1}(\mathbf w)\cap J_{3n,2}(\mathbf z)
\bigr|
\leq
|J_{2n,1}(\mathbf w)|
\leq
4b^{-(2+v_1)n}
\leq
4b^{-2N}.
\]
Thus the diameters of all intervals
in both coverings are at most $4b^{-2N}$,
which tends to zero as $N\to\infty$.

First, fix
\(s>\frac{2}{2+v_1}.
\)
Using the intervals $J_{2n,1}(\mathbf w)$
and \eqref{eq:exp-first-cover}, we obtain
\[
\begin{aligned}
\mathcal H^s\bigl(R_b(2,3,v_1,v_2)\bigr)
&\leq
\lim_{N\to\infty}
\sum_{n\geq N}
\ \sum_{\mathbf w\in\mathcal A_b^{2n}}
|J_{2n,1}(\mathbf w)|^s\\
&\leq
4^s\lim_{N\to\infty}
\sum_{n\geq N}
b^{[2-s(2+v_1)]n}
=0,
\end{aligned}
\]
because $2-s(2+v_1)<0$.
Since $s>2/(2+v_1)$ was arbitrary,
\[
\dimH R_b(2,3,v_1,v_2)
\leq\frac{2}{2+v_1}.
\]

Next, fix
\(s>
\frac{\max\{1,3-v_1,2-v_2\}}
{\max\{2+v_1,3+v_2\}}.
\)
By
\eqref{eq:exp-second-cover}
and 
\eqref{eq:exp-second-cover-cost}, we obtain
\[
\begin{aligned}
\mathcal H^s\bigl(R_b(2,3,v_1,v_2)\bigr)
&\quad\leq
\lim_{N\to\infty}
\sum_{n\geq N}
\ \sum_{\mathbf w\in\mathcal A_b^{2n}}
\ \sum_{\mathbf z\in\mathcal V_n(\mathbf w)}
\bigl|
J_{2n,1}(\mathbf w)\cap
J_{3n,2}(\mathbf z)
\bigr|^s\\
&\quad\leq
3C\,4^s\lim_{N\to\infty}
\sum_{n\geq N}
b^{n[
\max\{1,3-v_1,2-v_2\}
-s\max\{2+v_1,3+v_2\}
]}
=0.
\end{aligned}
\]
Since
\(s>
\frac{\max\{1,3-v_1,2-v_2\}}
{\max\{2+v_1,3+v_2\}}
\)
was arbitrary, we conclude that
\[
\dimH R_b(2,3,v_1,v_2)
\leq
\frac{\max\{1,3-v_1,2-v_2\}}
{\max\{2+v_1,3+v_2\}}.
\]
Combining the two upper bounds, we conclude
that
\begin{equation}\label{eq:23--up}
\dimH R_b(2,3,v_1,v_2)
\leq
\min\left\{
\frac{2}{2+v_1},
\frac{\max\{1,3-v_1,2-v_2\}}
{\max\{2+v_1,3+v_2\}}
\right\}.
\end{equation}

If $0\leq v_1\leq1$, then
\[
\max\{1,3-v_1,2-v_2\}=3-v_1,
\qquad
\max\{2+v_1,3+v_2\}=3+v_2.
\]
Substituting these identities into \eqref{eq:23--up}, we obtain
\[
\dimH R_b(2,3,v_1,v_2)
\leq \min\left\{
\frac{2}{2+v_1},
\frac{3-v_1}{3+v_2}
\right\}.
\]
If $v_1\geq1$, then
\[
\max\{1,3-v_1,2-v_2\}\leq2,
\qquad
\max\{2+v_1,3+v_2\}\geq2+v_1.
\]
Consequently,
\[
\frac{\max\{1,3-v_1,2-v_2\}}
{\max\{2+v_1,3+v_2\}}
\leq\frac{2}{2+v_1}.
\]
This proves the upper bound in Theorem~\ref{thm:exp}.

As an illustration, we point out that for $(v_1,v_2)=(2,1)$, the first cover gives
the bound $1/2$, whereas the second gives
the sharper bound $1/4$.

\section{The lower bound}\label{sec:exp-digits}
We prove the lower bound of Theorem~\ref{thm:exp} by constructing a Cantor type set
$C_{\infty}\subseteq R_b(2,3,v_1,v_2)$.
For this purpose, we choose a sequence of integers
$n_k\to\infty$ and construct digit sequences
satisfying
\[
|T_b^{2n_k}x-x|<b^{-v_1n_k},
\qquad
|T_b^{3n_k}x-x|<b^{-v_2n_k}
\]
for every $k\geq1$.
We ensure these inequalities by requiring the first $h \approx  v_1 n_k$ digits of $x$ and
\(T_b^{2n_k}x\) to agree, and the first $m\approx  v_2 n_k$ digits of $x$ and \(T_b^{3n_k}x\) to agree. 
These requirements are imposed
through the equalities of digits in the following
\eqref{eq:exp-copy}.
We first determine which positions in a
finite word can be chosen freely subject
to the equalities in \eqref{eq:exp-copy}.
Then we  use the resulting description to construct
the set $C_{\infty}$ and estimate its Hausdorff dimension.

\subsection{The free digit positions}

Let
\[
x=\sum_{i=1}^{\infty}w_i b^{-i},
\qquad w_i\in\mathcal A_b,
\]
be the  $b$-ary expansion of $x$.
Note that the $b$-ary exansion of $T_b^r x$ is
\((w_{r+1},w_{r+2},\ldots).
\)
We can impose
\(w_{2n+i}=w_i\) \((1\leq i\leq h)
\)
to make the first $h$ digits of $x$ and
$T_b^{2n}x$ agree. Similarly, we can impose
\(w_{3n+i}=w_i\) \((1\leq i\leq m)
\)
to make the first $m$ digits of $x$ and
$T_b^{3n}x$ agree. 
Then, we consider the two
families of digit equalities
\begin{equation}\label{eq:exp-copy}
w_{2n+i}=w_i\quad(1\leq i\leq h),
\qquad
w_{3n+i}=w_i\quad(1\leq i\leq m).
\end{equation}
The two families of equalities in
\eqref{eq:exp-copy} imply
\[
|T_b^{2n}x-x|<b^{-h},
\qquad
|T_b^{3n}x-x|<b^{-m}.
\]

The largest position appearing in \eqref{eq:exp-copy}
is
\[
L=\max\{2n+h,3n+m\}.
\]
Thus, the two families of equalities
in \eqref{eq:exp-copy} involve only the
first $L$ digits.
Therefore,
we  study finite words
$\mathbf w=(w_1,\ldots,w_L)\in\mathcal A_b^L$
satisfying \eqref{eq:exp-copy}.
We determine the positions at which the digits can
be prescribed freely, with all digits at the
other positions uniquely determined by
\eqref{eq:exp-copy}.
For each $1\leq\ell\leq L$, we also count
the freely prescribable positions among
the first $\ell$ digits.
\begin{lemma}\label{lem:exp-words}
Let $n\geq1$ and $h,m\geq0$ be integers, and put
\(L=\max\{2n+h,3n+m\}.
\)
When $h>n$, put
\(\kappa=\min\{h-n,m,n\}.
\)
Define
\begin{equation}\label{eq:exp-J}
\mathcal J=
\begin{cases}
\{1,\ldots,2n\}\cup
\{2n+h+1,\ldots,3n\},&h\leq n,\\[3pt]
\{1,\ldots,n\}\cup
\{n+\kappa+1,\ldots,2n\},&h>n.
\end{cases}
\end{equation}
Every assignment of digits from
$\mathcal A_b=\{0,\ldots,b-1\}$ to the positions
in $\mathcal J$ determines a unique word
$\mathbf w=(w_1,\ldots,w_L)$ satisfying
\eqref{eq:exp-copy}.
For each integer $1\leq\ell\leq L$, define
\[
f(\ell):=\#\bigl(\mathcal J\cap\{1,\ldots,\ell\}\bigr).
\]
Then exactly $b^{f(\ell)}$ distinct prefixes
of length $\ell$ occur among the words
satisfying \eqref{eq:exp-copy}

Finally, with
\[
v_1=\frac hn,\qquad v_2=\frac mn,
\]
the minimum proportion of freely prescribable
positions is
\begin{equation}\label{eq:exp-min}
\min_{1\leq\ell\leq L}\frac{f(\ell)}{\ell}
=
\min_{n\leq\ell\leq L}\frac{f(\ell)}{\ell}
=D_{2,3}(v_1,v_2).
\end{equation}
\end{lemma}

\begin{figure}[htbp]
\centering
\begin{tikzpicture}[
  x=0.85cm,
  y=1cm,
  line cap=round,
  line join=round,
  every node/.style={font=\small, inner sep=2pt},
  free/.style={draw=black, fill=white, line width=.55pt},
  copied/.style={draw=black, fill=gray!25, line width=.55pt}
]

\def\hh{0.46}

\node[anchor=west, font=\normalsize] at (0, 1.15) {(a) $h\leq n$};

\draw[free]   (0,   -\hh) rectangle (3.4,  \hh);
\draw[copied] (3.4, -\hh) rectangle (6.9,  \hh);
\draw[free]   (6.9, -\hh) rectangle (10.2, \hh);
\draw[copied] (10.2,-\hh) rectangle (14.0, \hh);

\node[align=center, font=\footnotesize] at (1.7,  0) {freely\\prescribable};
\node                                   at (5.15, 0) {$w_i=w_{i-2n}$};
\node[align=center, font=\footnotesize] at (8.55, 0) {freely\\prescribable};
\node                                   at (12.1, 0) {$w_i=w_{i-3n}$};

\foreach \x/\lab in {
  0/0,
  3.4/{2n},
  6.9/{2n+h},
  10.2/{3n},
  14.0/L
}{
  \draw[line width=.55pt] (\x, -\hh) -- ++(0, -0.12);
  \node[below=3pt] at (\x, -\hh-0.12) {$\lab$};
}

\node[anchor=west, font=\normalsize] at (0, -1.65) {(b) $h>n$};

\def\yB{-2.85}

\draw[free]   (0,   \yB-\hh) rectangle (2.5,  \yB+\hh);
\draw[copied] (2.5, \yB-\hh) rectangle (5.1,  \yB+\hh);
\draw[free]   (5.1, \yB-\hh) rectangle (7.7,  \yB+\hh);
\draw[copied] (7.7, \yB-\hh) rectangle (10.6, \yB+\hh);
\draw[copied] (10.6,\yB-\hh) rectangle (14.0, \yB+\hh);

\node[align=center, font=\footnotesize] at (1.25, \yB) {freely\\prescribable};
\node                                   at (3.8,  \yB) {$w_i=w_{i-n}$};
\node[align=center, font=\footnotesize] at (6.4,  \yB) {freely\\prescribable};
\node                                   at (9.15, \yB) {$w_i=w_{i-2n}$};
\node[align=center, font=\footnotesize] at (12.3, \yB) {$w_i=w_{i-2n}$\\or $w_{i-3n}$};

\foreach \x/\lab in {
  0/0,
  2.5/n,
  5.1/{n+\kappa},
  7.7/{2n},
  10.6/{3n+t},
  14.0/L
}{
  \draw[line width=.55pt] (\x, \yB-\hh) -- ++(0, -0.12);
  \node[below=3pt] at (\x, \yB-\hh-0.12) {$\lab$};
}

\end{tikzpicture}
\caption{Freely prescribable and determined positions in
Lemma~\ref{lem:exp-words} for (a) $h\leq n$ and (b) $h>n$, where
$t=\min\{h-n,m\}$. White blocks represent the positions in $\mathcal J$,
and shaded blocks represent positions determined by \eqref{eq:exp-copy}
through the indicated equalities. The widths are schematic, and some
blocks may be empty.}
\label{fig:exp-free-positions}
\end{figure}
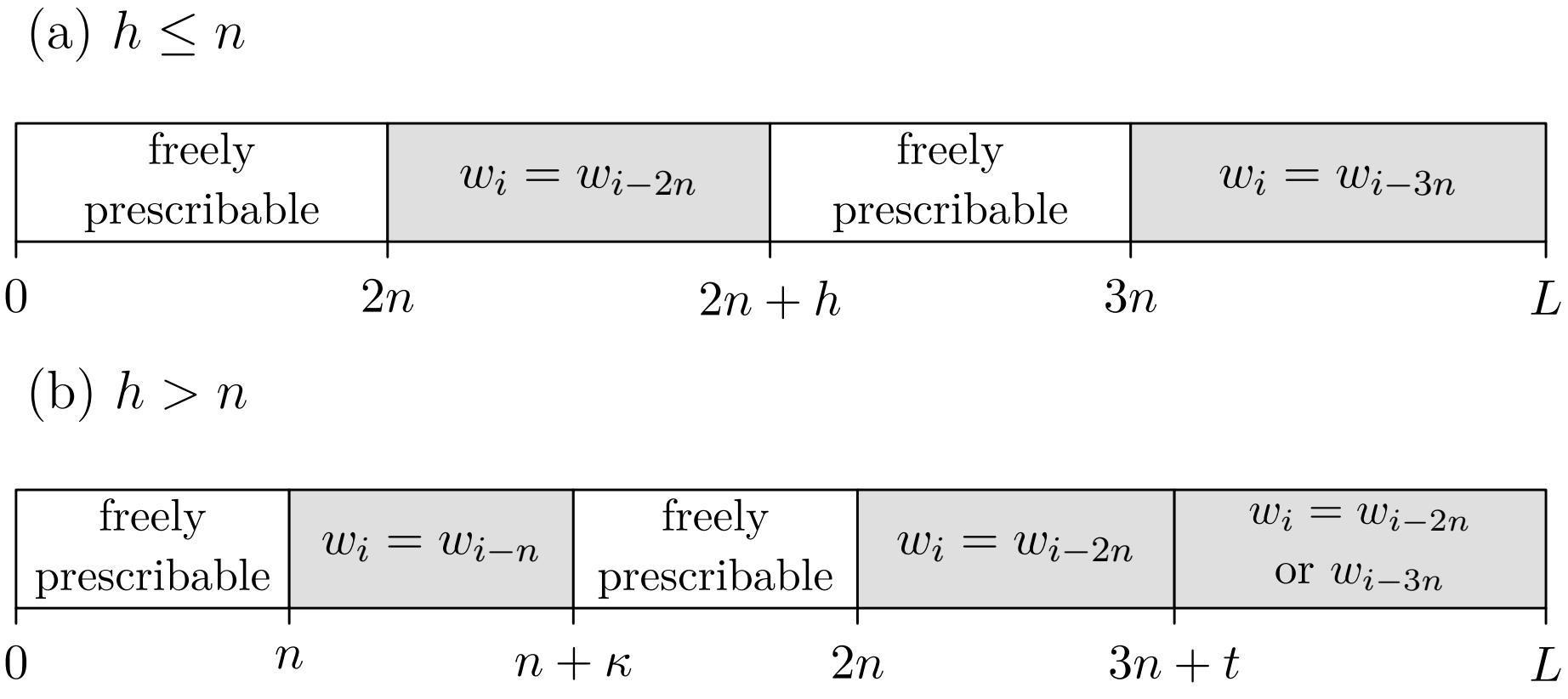

Figure~\ref{fig:exp-free-positions}
illustrates the two forms of $\mathcal J$
defined in \eqref{eq:exp-J}.

\begin{proof}
For each $i\in\mathcal J$, choose a digit
$w_i\in\mathcal A_b$ arbitrarily.
We show that the remaining digits
$w_i$, $i\in\{1,\ldots,L\}\setminus\mathcal J$,
can be assigned in exactly one way so that
$\mathbf w=(w_1,\ldots,w_L)$ satisfies
\eqref{eq:exp-copy}.

\medskip
\noindent
\emph{Case 1: $h\leq n$.}
Since $2n+h\leq3n$,
the two sets of positions
\[
\{2n+1,\ldots,2n+h\}
\quad\text{and}\quad
\{3n+1,\ldots,3n+m\}
\]
are disjoint.
By \eqref{eq:exp-J}, the digits at positions
\[
\{1,\ldots,2n\}
\cup\{2n+h+1,\ldots,3n\}
\]
have already been prescribed.
For $2n+1\leq i\leq2n+h$, define
\(w_i=w_{i-2n}.
\)
Since $1\leq i-2n\leq h\leq n$,
	the digit $w_{i-2n}$ has already been
	chosen.	
Now, all digits at positions $1,\ldots,3n$
	are  defined.
	For $j=1,\ldots,m$, in increasing order,
	define
	\[
	w_{3n+j}=w_j.
	\]
	At each step, $j<3n+j$, so the digit
	$w_j$ has already been defined.
Thus all digits outside $\mathcal J$
have been uniquely determined, and
both families of equalities in
\eqref{eq:exp-copy} are satisfied.

\medskip
\noindent
\emph{Case 2: $h>n$.}
Put
\[
t=\min\{h-n,m\},
\qquad
\kappa=\min\{t,n\}.
\]
For $1\leq i\leq t$, we have
$n+i\leq h$ and $i\leq m$.
Thus the two families of equalities
in \eqref{eq:exp-copy} give
\[
w_{3n+i}=w_{n+i}
\quad\text{and}\quad
w_{3n+i}=w_i.
\]
Consequently, every word satisfying
\eqref{eq:exp-copy} must satisfy
\begin{equation}\label{eq:exp-overlap-compatibility}
w_{n+i}=w_i
\qquad(1\leq i\leq t).
\end{equation}

By \eqref{eq:exp-J}, the digits at positions
\[
\mathcal J
=\{1,\ldots,n\}
\cup\{n+\kappa+1,\ldots,2n\}
\]
have already been prescribed.
We first define, in increasing order of $i$,
\[
w_{n+i}=w_i,
\qquad 1\leq i\leq t.
\]
If $i\leq n$, the digit $w_i$ was prescribed
at the beginning of the proof.
If $n<i\leq t$, the digit $w_i$
was determined at the earlier step $i-n$.
Since $\kappa=\min\{t,n\}$,
the positions $n+1,\ldots,n+t$
are disjoint from $\mathcal J$:
if $t\leq n$, then $\kappa=t$;
if $t>n$, then $\mathcal J=\{1,\ldots,n\}$.
Thus the assignments $w_{n+i}=w_i$
do not change any of the digits
prescribed on $\mathcal J$.

If $t\leq n$, then $\kappa=t$ and
\(\mathcal J
=\{1,\ldots,n\}\cup\{n+t+1,\ldots,2n\}.
\)
The digits at positions in $\mathcal J$
were prescribed at the beginning of the proof,
while the equalities $w_{n+i}=w_i$ for
$1\leq i\leq t$ determine the digits
at positions $n+1,\ldots,n+t$.
Hence all digits $w_1,\ldots,w_{2n}$
have been assigned.
If $t>n$, then $\kappa=n$ and
$\mathcal J=\{1,\ldots,n\}$.
The digits $w_1,\ldots,w_n$ were prescribed,
and the equalities $w_{n+i}=w_i$ for
$1\leq i\leq t$ determine the digits
$w_{n+1},\ldots,w_{n+t}$.
Hence all digits $w_1,\ldots,w_{n+t}$
have been assigned.
Thus, putting
\[
M=\max\{2n,n+t\},
\]
all digits $w_1,\ldots,w_M$ have been assigned.

We next impose the first family of equalities
in \eqref{eq:exp-copy} for
$1\leq i\leq n+t$.
For $1\leq i\leq t-n$, when $t>n$,
the digit $w_{2n+i}$ has already been
assigned, and
\[
w_{2n+i}=w_{n+i}=w_i
\]
follows from
\eqref{eq:exp-overlap-compatibility}.
For
$\max\{1,t-n+1\}\leq i\leq n+t$,
we have
\[
2n+i
\geq 2n+\max\{1,t-n+1\}
=M+1,
\]
and $i\leq n+t\leq M$.
Thus $w_i$ has already been assigned,
whereas $w_{2n+i}$ has not.

All digits $w_1,\ldots,w_{3n+t}$
have now been assigned, and
\[
w_{2n+i}=w_i
\qquad(1\leq i\leq n+t).
\]
Furthermore, for $1\leq i\leq t$,
\[
w_{3n+i}=w_{n+i}=w_i.
\]
The first equality holds because
$n+i\leq n+t$, and the second follows
from \eqref{eq:exp-overlap-compatibility}.
Hence the second family of equalities
in \eqref{eq:exp-copy} is also satisfied
for $1\leq i\leq t$.

It remains to assign the digits at positions
greater than $3n+t$.

If $m<h-n$, then $t=m$ and $L=2n+h$.
For $3n+t<i\leq L$, define
\(w_i=w_{i-2n}
\)
in increasing order of $i$.
This imposes the remaining equalities
in the first family of \eqref{eq:exp-copy}.
If $h-n<m$, then $t=h-n$ and $L=3n+m$.
For $3n+t<i\leq L$, define
\[
w_i=w_{i-3n}
\]
in increasing order of $i$.
This imposes the remaining equalities
in the second family of \eqref{eq:exp-copy}.
If $m=h-n$, then $L=3n+t$,
so no further assignments are needed.
In the first two cases, each newly assigned
digit is copied from an earlier position,
so assigning the digits in increasing
order of their indices is well-defined.

The constructions in Cases 1 and 2 show that
every assignment of digits from $\mathcal A_b$
to the positions in $\mathcal J$ determines
exactly one word of length $L$ satisfying
\eqref{eq:exp-copy}.
In both constructions, every digit $w_i$
with $i\notin\mathcal J$ is determined by
an equality $w_i=w_j$ with $1\leq j<i$.
The index $j$ depends only on $n,h,m$ and $i$,
not on the values of the prescribed digits.
It follows by induction on $i$ that $w_i$
is determined by the digits at positions in
$\mathcal J\cap\{1,\ldots,i\}$.
Therefore, for every integer $1\leq\ell\leq L$,
the digits at positions in
\[
\mathcal J\cap\{1,\ldots,\ell\}
\]
uniquely determine the prefix
$(w_1,\ldots,w_\ell)$.

Conversely, fix any assignment of digits
to the positions in
$\mathcal J\cap\{1,\ldots,\ell\}$.
Assign arbitrary digits to the positions in
$\mathcal J\cap\{\ell+1,\ldots,L\}$.
The constructions in Cases 1 and 2
then uniquely determine the digits outside
$\mathcal J$ and yield a word satisfying
\eqref{eq:exp-copy}.
Any two distinct assignments of digits to
$\mathcal J\cap\{1,\ldots,\ell\}$
differ at some position
$i\in\mathcal J\cap\{1,\ldots,\ell\}$,
and therefore yield distinct prefixes
of length $\ell$.
Hence the number of distinct prefixes
of length $\ell$ is exactly
\[
b^{\#(\mathcal J\cap\{1,\ldots,\ell\})}
=b^{f(\ell)}.
\]

We now determine the minimum of
$f(\ell)/\ell$ over integers
$1\leq\ell\leq L$.
Since $0\leq f(\ell)\leq\ell$, for
$1\leq\ell<L$, if $\ell+1\in \mathcal J$, then
\[
\frac{f(\ell+1)}{\ell+1}
=\frac{f(\ell)+1}{\ell+1}
\geq\frac{f(\ell)}{\ell}.
\]
If $\ell+1\notin \mathcal J$, then
\[
\frac{f(\ell+1)}{\ell+1}
=\frac{f(\ell)}{\ell+1}
\leq\frac{f(\ell)}{\ell}.
\]

Suppose first that $h\leq n$.
Then $L=3n+m$.
By \eqref{eq:exp-J}, the freely prescribable positions
are
\[
\mathcal J=
\{1,\ldots,2n\}
\cup
\{2n+h+1,\ldots,3n\}.
\]
The positions $2n+1,\ldots,2n+h$ are not
freely prescribable. 
The ratio $f(\ell)/\ell$ equals $1$ for
$1\leq\ell\leq2n$, is nonincreasing for
$2n\leq\ell\leq2n+h$, nondecreasing for
$2n+h\leq\ell\leq3n$, and nonincreasing
for $3n\leq\ell\leq L$.
Thus the minimum is attained at
$\ell=2n+h$ or $\ell=L$.
Since
\(f(2n+h)=2n\) and
\(f(L)=3n-h
\),
we obtain
\[
\begin{aligned}
\min_{1\leq\ell\leq L}\frac{f(\ell)}{\ell}
&=
\min\left\{
\frac{2n}{2n+h},
\frac{3n-h}{3n+m}
\right\}\\
&=
\min\left\{
\frac{2}{2+v_1},
\frac{3-v_1}{3+v_2}
\right\}
=D_{2,3}(v_1,v_2),
\end{aligned}
\]
where $v_1=h/n$ and $v_2=m/n$.

Now, suppose  that $h>n$.
By \eqref{eq:exp-J},
\[
\mathcal J
=\{1,\ldots,n\}
\cup\{n+\kappa+1,\ldots,2n\}.
\]
The ratio $f(\ell)/\ell$ equals $1$ for
$1\leq\ell\leq n$, is nonincreasing for
$n\leq\ell\leq n+\kappa$, nondecreasing
for $n+\kappa\leq\ell\leq2n$, and
nonincreasing for $2n\leq\ell\leq L$.
Since
\[
f(n+\kappa)=n,
\qquad
f(L)=2n-\kappa,
\]
we obtain
\[
\min_{1\leq\ell\leq L}\frac{f(\ell)}{\ell}
=
\min\left\{
\frac{n}{n+\kappa},
\frac{2n-\kappa}{L}
\right\}.
\]

Since $0\leq\kappa\leq n$ and $L\geq3n$,
\[
\begin{aligned}
(2n-\kappa)(n+\kappa)
=2n^2+n\kappa-\kappa^2\leq3n^2
\leq nL.
\end{aligned}
\]
Dividing by $L(n+\kappa)>0$ gives
\(\frac{2n-\kappa}{L}
\leq\frac{n}{n+\kappa}.
\)
Thus the minimum is attained at $\ell=L$.

By the definitions of $\kappa,L,v_1,v_2$,
\[
\frac{\kappa}{n}
=\min\{v_1-1,v_2,1\},
\qquad
\frac Ln=\max\{2+v_1,3+v_2\}.
\]
Consequently,
\[
\begin{aligned}
\min_{1\leq\ell\leq L}\frac{f(\ell)}{\ell}
&=\frac{2n-\kappa}{L}\\
&=
\frac{\max\{1,3-v_1,2-v_2\}}
{\max\{2+v_1,3+v_2\}}\\
&=D_{2,3}(v_1,v_2),
\end{aligned}
\]
where the last equality uses $v_1>1$.

When $h\leq n$, the minimum is attained
at $2n+h$ or $L$; when $h>n$, the minimum
is attained at $L$.
Since $2n+h\geq n$ and $L\geq3n$,
restricting $\ell$ to $\{n,\ldots,L\}$
does not change the minimum.
This proves \eqref{eq:exp-min}.
\end{proof}

The minimum in \eqref{eq:exp-min} need not
be attained at $\ell=L$.
For example, take $h=n/2$ and $m=n/16$,
where $n$ is a positive multiple of $16$.
Then $L=49n/16$, and \eqref{eq:exp-J} gives
\[
\mathcal J
=
\{1,\ldots,2n\}
\cup
\{5n/2+1,\ldots,3n\}.
\]
At the two lengths $2n+h$ and $L$, we have
\[
\begin{aligned}
\frac{f(2n+h)}{2n+h}
=\frac{2n}{5n/2}
=\frac45,\quad
\frac{f(L)}{L}
=\frac{5n/2}{49n/16}
=\frac{40}{49}.
\end{aligned}
\]
Since $4/5<40/49$, the proof of
Lemma~\ref{lem:exp-words} shows that
the minimum in \eqref{eq:exp-min}
is attained at
\[
\ell=2n+h=\frac{5n}{2}
<
L=\frac{49n}{16}.
\]

Although $\mathcal J$ is defined in terms of
$h$ and $m$, different pairs $(h,m)$ may
give the same set $\mathcal J$.
The following example shows that such pairs
can produce different words even when the
digits prescribed on $\mathcal J$ are identical.
For example, let $n=2$ and
	choose $a_1,a_2,a_3\in\mathcal A_b$ with
	$a_2\ne a_3$. Prescribe
	$w_1=a_1$, $w_2=a_2$, and $w_4=a_3$.
Both choices
\[
(h,m)=(3,2)
\qquad\text{and}\qquad
(h,m)=(4,1)
\]
give the same set
\(\mathcal J=\{1,2,4\}.
\)
	For both choices, $L=8$, and
	\eqref{eq:exp-copy} gives
\[
(w_1,\ldots,w_7)
=(a_1,a_2,a_1,a_3,a_1,a_2,a_1).
\]

For $(h,m)=(3,2)$, the second family of equalities in
\eqref{eq:exp-copy}, applied with $i=2$,
gives
\[
w_8=w_2=a_2.
\]
For $(h,m)=(4,1)$, the first family of equalities in
\eqref{eq:exp-copy}, applied with $i=4$,
gives
\[
w_8=w_4=a_3.
\]
Since $a_2\ne a_3$, the words corresponding to
$(h,m)=(3,2)$ and $(h,m)=(4,1)$ differ
at position $8$, although both have the same
set of freely prescribable positions
$\mathcal J=\{1,2,4\}$ and the same prescribed
digits $w_1=a_1$, $w_2=a_2$, and $w_4=a_3$.

\subsection{Construction of the subset}
\label{subsec:exp-lower}

We first assume that $v_1$ and $v_2$ are positive
rational numbers and put
\[
D=D_{2,3}(v_1,v_2).
\]
We construct finite families of words
$\mathcal F_k\subseteq\mathcal A_b^{L_k}$,
$k\geq0$, such that the prefix of length
$L_{k-1}$ of every word in $\mathcal F_k$
belongs to $\mathcal F_{k-1}$.

We start with
\[
L_0=1,\qquad \mathcal F_0=\{(0)\}.
\]

Suppose that $k\geq1$ and that
$L_{k-1}$ and $\mathcal F_{k-1}$ have been defined.
Choose $n_k$ to be a sufficiently large
positive multiple of a common denominator
of $v_1$ and $v_2$ so that
\[
n_k>k^2L_{k-1}.
\]
Set
\begin{equation}\label{expd:times}
h_k=v_1n_k,\qquad
m_k=v_2n_k,\qquad
L_k=\max\{2n_k+h_k,3n_k+m_k\}.
\end{equation}
Then $h_k$ and $m_k$ are positive integers, and
\[
L_{k-1}<n_k<L_k,
\qquad
\frac{L_{k-1}}{n_k}<\frac1{k^2}.
\]

At stage $k$, we require the digits of each
word $\mathbf w=(w_1,\ldots,w_{L_k})$ to satisfy
\begin{equation}\label{expd:copies}
w_{2n_k+i}=w_i\quad(1\leq i\leq h_k),
\qquad
w_{3n_k+i}=w_i\quad(1\leq i\leq m_k).
\end{equation}
By Lemma~\ref{lem:exp-words}, the positions
at which digits can be prescribed freely are
\begin{equation}\label{expd:positions}
\mathcal J_k=
\begin{cases}
\{1,\ldots,2n_k\}
\cup\{2n_k+h_k+1,\ldots,3n_k\},
&h_k\leq n_k,\\[3pt]
\{1,\ldots,n_k\}
\cup\{n_k+\kappa_k+1,\ldots,2n_k\},
&h_k>n_k,
\end{cases}
\end{equation}
where
\[
\kappa_k=\min\{h_k-n_k,m_k,n_k\}
\qquad\text{if }h_k>n_k.
\]
An indicated set is understood to be empty
when its lower endpoint exceeds its upper endpoint.

For $1\leq\ell\leq L_k$, let
\[
f_k(\ell)
=\#\bigl(\mathcal J_k\cap\{1,\ldots,\ell\}\bigr).
\]
Since
$\{1,\ldots,n_k\}\subseteq \mathcal J_k$,
Lemma~\ref{lem:exp-words} gives
\begin{equation}\label{eq:exp-free-density}
f_k(\ell)=\ell\quad(1\leq\ell\leq n_k),
\qquad
f_k(\ell)\geq D\ell\quad(1\leq\ell\leq L_k).
\end{equation}

We now define the words in $\mathcal F_k$.
Fix a word
\(W=(W_1,\ldots,W_{L_{k-1}})
\in\mathcal F_{k-1}.
\)
Since $L_{k-1}<n_k$ and
$\{1,\ldots,n_k\}\subseteq \mathcal J_k$,
	all positions $1,\ldots,L_{k-1}$ belong
	to $\mathcal J_k$.
	We prescribe the prefix of length
	$L_{k-1}$ of the new word by setting
	\[
	w_i=W_i\qquad(1\leq i\leq L_{k-1}).
	\]
We then assign arbitrary digits from
	$\mathcal A_b$ to the remaining positions
	in $\mathcal J_k$, namely,
	\[
	\mathcal J_k\cap
	\{L_{k-1}+1,\ldots,L_k\}.
	\]
	By Lemma~\ref{lem:exp-words}, these
	assignments determine all digits outside
	$\mathcal J_k$ uniquely.

For each $W\in\mathcal F_{k-1}$, define
\begin{equation}\label{expd:children}
\mathcal F_k(W)=
\left\{
\mathbf w=(w_1,\ldots,w_{L_k}):
\begin{array}{ll}
w_i=W_i,&1\leq i\leq L_{k-1},\\
w_i\in\mathcal{A}_b,
&i\in J_k\cap\{L_{k-1}+1,\ldots,L_k\},\\
\mathbf w\text{ satisfies }\eqref{expd:copies}.
\end{array}
\right\}.
\end{equation}
Set
\[
\mathcal F_k
=\bigcup_{W\in\mathcal F_{k-1}}\mathcal F_k(W).
\]
Every word in $\mathcal F_k(W)$ has
$W$ as its prefix of length $L_{k-1}$.
Since distinct words in $\mathcal F_{k-1}$
have distinct prefixes of length $L_{k-1}$,
every word in $\mathcal F_k$ has a unique
parent word in $\mathcal F_{k-1}$.

Since $L_{k-1}<n_k$ and
$\{1,\ldots,n_k\}\subseteq\mathcal J_k$,
we have
$f_k(L_{k-1})=L_{k-1}$.
For a fixed $W\in\mathcal F_{k-1}$,
the number of positions in $\mathcal J_k$
at which new digits can be chosen is
\[
f_k(L_k)-f_k(L_{k-1})
=f_k(L_k)-L_{k-1}.
\]
Each position in
$\mathcal J_k\cap
\{L_{k-1}+1,\ldots,L_k\}$
admits exactly $b$ choices.
By Lemma~\ref{lem:exp-words}, every
assignment of digits to these positions,
together with the fixed prefix $W$,
determines a unique word in
$\mathcal F_k(W)$.
Consequently,
\begin{equation}\label{expd:branching}
\#\mathcal F_k(W)
=b^{f_k(L_k)-L_{k-1}}
\qquad(W\in\mathcal F_{k-1}).
\end{equation}

To illustrate how the prefix condition
is preserved from one stage to the next,
consider $v_1=2$ and $v_2=1$.
Then
\(h_k=2n_k\), \( m_k=n_k,
\), \(\mathcal J_k=\{1,\ldots,n_k\}\) and 
\(L_k=4n_k.
\)
At the first stage, choose a word $S_1$
of length $n_1$ whose first digit is $0$.
The two digit equalities \eqref{eq:exp-copy} then give
\[
W_1=S_1^4.
\]
At the second stage, choose a word $S_2$
of length $n_2$ whose first $4n_1$ digits
agree with all $4n_1$ digits of $W_1$.
This is possible because
$n_2>4L_1=16n_1$.
The remaining $n_2-4n_1$ digits of $S_2$
can be chosen arbitrarily. The equalities
in \eqref{expd:copies} then determine
\[
W_2=S_2^4.
\]
Thus the first $4n_1$ digits of $W_2$
are exactly the word $W_1$.

 For each $k\geq0$,
put
\begin{equation}\label{expd:limit}
C_k=
\bigcup_{W\in\mathcal F_k}
I_{L_k}(W),
\qquad
C_{\infty}=\bigcap_{k\geq0}C_k.
\end{equation}
In particular,
\(C_0=I_{1}(0)=[0,b^{-1}).
\)

For each $W\in\mathcal F_{k-1}$,
the intervals
\(I_{L_k}(\mathbf w)\),
 \(\mathbf w\in\mathcal F_k(W),
\)
are contained in $I_{L_{k-1}}(W)$,
because every word $\mathbf w\in\mathcal F_k(W)$
has $W$ as its prefix of length $L_{k-1}$.
There are exactly
$b^{f_k(L_k)-L_{k-1}}$
such intervals, and each interval has
length $b^{-L_k}$.
Taking the union over
$W\in\mathcal F_{k-1}$ gives
\[
C_k\subseteq C_{k-1}
\qquad(k\geq1).
\]

Finally, let $x\in C_{\infty}$ and write its 
$b$-expansion as
\(x=\sum_{i=1}^{\infty}w_i b^{-i}.
\)
Since $C_k\subseteq C_{k-1}$ for every
$k\geq1$, and the cylinders of order
$L_k$ are pairwise disjoint, each
$x\in C_\infty$ belongs to a unique
cylinder $I_{L_k}(W_k)$ with
$W_k\in\mathcal F_k$.
By the definition of $\mathcal F_k$, the digits
of $x$ satisfy
\[
w_{2n_k+i}=w_i\qquad(1\leq i\leq h_k),
\]
and
\[
w_{3n_k+i}=w_i\qquad(1\leq i\leq m_k).
\]
Therefore, the first $h_k$ digits of $x$
and $T_b^{2n_k}x$ agree, while the first
$m_k$ digits of $x$ and $T_b^{3n_k}x$ agree.
It follows that
\begin{equation}\label{expd:returns}
\begin{aligned}
|T_b^{2n_k}x-x|
<b^{-h_k}=b^{-v_1n_k},\quad
|T_b^{3n_k}x-x|
<b^{-m_k}=b^{-v_2n_k}.
\end{aligned}
\end{equation}
Since $n_k\to\infty$, both inequalities hold
simultaneously for infinitely many integers $n_k$.
Consequently,
\[
C_{\infty}\subseteq R_b(2,3,v_1,v_2).
\]

\subsection{The measure on the constructed set}
\label{subsec:exp-measure}

\smallskip
\noindent
\emph{Definition of the measure}
We assign mass one to the initial interval
$I_1(0)$ and mass zero to every other
cylinder of order one.

For each $k\geq1$, every word
$W\in\mathcal F_{k-1}$ has exactly
$b^{f_k(L_k)-L_{k-1}}$ children in
$\mathcal F_k(W)$.
We assign equal mass to the cylinders of order $L_k$
corresponding to these children.
More precisely,
for $W_k\in\mathcal F_k(W)$, then
\begin{equation}\label{eq:expd-child-mass}
\mu(I_{L_k}(W_k))
=
b^{-[f_k(L_k)-L_{k-1}]}
\mu(I_{L_{k-1}}(W)).
\end{equation}

For every integer $\ell\geq2$, there is
a unique integer $k\geq1$ such that
$L_{k-1}<\ell\leq L_k$.
For each cylinder $I_\ell(\mathbf u)$
of order $\ell$ that meets $C_\infty$,
define
\begin{equation}\label{eq:expd-general-mass}
\mu(I_{\ell}(\mathbf u))
=
\sum_{\substack{W_k\in\mathcal F_k\\
I_{L_k}(W_k)\subseteq I_{\ell}(\mathbf u)}}
\mu(I_{L_k}(W_k)).
\end{equation}
Cylinders that do not meet $C_\infty$
are assigned mass zero.
The mass of each cylinder of order $\ell$
equals the sum of the masses of the $b$
cylinders of order $\ell+1$ contained
in that cylinder.
Indeed, the child masses in
\eqref{eq:expd-child-mass} sum to the
mass of their parent cylinder.
Thus the cylinder masses are consistent
across all orders.
By the standard extension theorem for
consistent masses on nested cylinder
partitions, they determine a unique
Borel probability measure $\mu$ on $[0,1)$.
For every $k\geq0$, the total mass of
the cylinders $I_{L_k}(W)$,
$W\in\mathcal F_k$, is one.
Hence $\mu(C_k)=1$.
Since $C_k$ is decreasing, continuity
of measures from above gives
\[
\mu(C_\infty)
=\mu\left(\bigcap_{k\geq0}C_k\right)
=\lim_{k\to\infty}\mu(C_k)=1.
\]

\smallskip
\noindent
\emph{Estimation of the masses of intervals}
Fix $k\geq1$ and
$L_{k-1}<\ell\leq L_k$.
Let $x\in C_\infty$, and let
$W\in\mathcal F_{k-1}$ be the unique word
such that
$x\in I_{L_{k-1}}(W)$.

Among the positions $L_{k-1}+1,\ldots,\ell$,
exactly $f_k(\ell)-L_{k-1}$ positions
belong to $\mathcal J_k$.
The digits at these positions are fixed
by the prefix of length $\ell$.
The remaining freely prescribable positions
at stage $k$ are precisely
\(\mathcal J_k\cap\{\ell+1,\ldots,L_k\},
\)
whose number is $f_k(L_k)-f_k(\ell)$.
Consequently, exactly
\[
b^{f_k(L_k)-f_k(\ell)}
\]
children of $W$ have the same first $\ell$
digits as $x$. Since all children of $W$
have equal mass, \eqref{eq:expd-child-mass}
gives
\begin{equation}\label{eq:expd-prefix-mass}
\mu(I_\ell(x))
=
\mu(I_{L_{k-1}}(x))
b^{-[f_k(\ell)-L_{k-1}]}.
\end{equation}

For each integer $\ell\geq1$, let $G(\ell)$
denote the total number of freely chosen
digits among the first $\ell$ positions
through all stages of the construction.
The first digit is fixed as $0$ by
$\mathcal F_0=\{(0)\}$, so $G(1)=0$.
By \eqref{eq:expd-prefix-mass},
\[
G(\ell)=G(L_{k-1})+
f_k(\ell)-L_{k-1},
\]
and for any $x\in C_{\infty}$,
\begin{equation}\label{eq:exp-cumulative}
\mu(I_\ell(x))=b^{-G(\ell)}.
\end{equation}

We show that
\[
\liminf_{\ell\to\infty}
\frac{G(\ell)}{\ell}\geq D.
\]

First, suppose that $n_k\leq\ell\leq L_k$.
By \eqref{eq:exp-free-density},
$f_k(\ell)\geq D\ell$.
Since $G(L_{k-1})\geq0$, we obtain
\[
\begin{aligned}
\frac{G(\ell)}{\ell}
&=\frac{G(L_{k-1})+f_k(\ell)-L_{k-1}}{\ell}\\
&\geq D-\frac{L_{k-1}}{\ell}\\
&\geq D-\frac{L_{k-1}}{n_k}.
\end{aligned}
\]

Next, suppose that $L_{k-1}\leq\ell\leq n_k$.
Since
$\{1,\ldots,n_k\}\subseteq\mathcal J_k$
and the first $L_{k-1}$ digits were fixed
by the parent word, all positions
$L_{k-1}+1,\ldots,\ell$ are freely
assigned at stage $k$.
 Thus
\[
G(\ell)
=G(L_{k-1})+\ell-L_{k-1}.
\]
It follows that
\[
\begin{aligned}
\frac{G(\ell)}{\ell}
&=1-\frac{L_{k-1}-G(L_{k-1})}{\ell}\\
&\geq
1-\frac{L_{k-1}-G(L_{k-1})}{L_{k-1}}\\
&=\frac{G(L_{k-1})}{L_{k-1}}.
\end{aligned}
\]
The inequality follows from
$G(L_{k-1})\leq L_{k-1}$ and
$\ell\geq L_{k-1}$.

Applying the estimate for the interval
$n_{k-1}\leq\ell\leq L_{k-1}$
at $\ell=L_{k-1}$ gives
\[
\frac{G(L_{k-1})}{L_{k-1}}
\geq D-\frac{L_{k-2}}{n_{k-1}}.
\]
Combining the estimates for
$L_{k-1}\leq\ell\leq n_k$
and $n_k\leq\ell\leq L_k$,
we obtain
\[
\frac{G(\ell)}{\ell}
\geq
D-\max\left\{
\frac{L_{k-1}}{n_k},
\frac{L_{k-2}}{n_{k-1}}
\right\}
\]
for $L_{k-1}\leq\ell\leq L_k$.
By the choice of the integers $n_k$,
\[
\frac{L_{k-1}}{n_k}<\frac1{k^2},
\qquad
\frac{L_{k-2}}{n_{k-1}}<\frac1{(k-1)^2}.
\]
Since $L_k>L_{k-1}$ for every $k\geq1$
and $L_k\to\infty$, the intervals
$[L_{k-1},L_k]$ cover all integers
$\ell\geq L_0$.
Therefore,
\[
\liminf_{\ell\to\infty}
\frac{G(\ell)}{\ell}\geq D.
\]

Fix $s$ with $0<s<D$.
Since
$\liminf_{\ell\to\infty}G(\ell)/\ell\geq D$,
there exists $\ell_0$ such that for any $\ell\geq\ell_0$,
\[
G(\ell)\geq s\ell.
\]
It follows from \eqref{eq:exp-cumulative} that
\begin{equation}\label{eq:exp-cylinder-estimate}
\mu(I_\ell(x))\leq b^{-s\ell}
\qquad(x\in C_{\infty},\ \ell\geq\ell_0).
\end{equation}

\smallskip
\noindent
\emph{Estimation for arbitrary intervals.}
We first note that $\mu$ has no atoms.
Indeed, for every $x\in C_\infty$,
\eqref{eq:exp-cylinder-estimate} gives for any $\ell\geq\ell_0$,
\[
\mu(\{x\})
\leq\mu(I_\ell(x))
\leq b^{-s\ell}.
\]
Letting $\ell\to\infty$ yields
$\mu(\{x\})=0$.
Since $\mu(C_\infty)=1$, every point
outside $C_\infty$ also has zero mass.

Let $I\subset\mathbb R$ be an interval
of sufficiently small positive length.
Choose an integer $\ell\geq\ell_0$ such that
\[
b^{-(\ell+1)}<|I|\leq b^{-\ell}.
\]
The cylinders of order $\ell$ form a
partition of $[0,1)$ into pairwise disjoint
intervals of length $b^{-\ell}$.
Since $|I|\leq b^{-\ell}$,
$I$ intersects at most two cylinders
in sets of positive length.
Because $\mu$ has no atoms, intersections
consisting only of endpoints contribute
zero mass.
For the remaining cylinders,
\eqref{eq:exp-cylinder-estimate} gives
\[
\mu(I)\leq2b^{-s\ell}
\leq2b^s|I|^s.
\]
By Proposition~\ref{mmd}, the estimate
$\mu(I)\leq2b^s|I|^s$ for all sufficiently
small intervals $I$ gives
\[
\dim_H C_\infty\geq s.
\]
Since $C_\infty\subseteq R_b(2,3,v_1,v_2)$,
we obtain
\(\dim_{\mathrm{H}} R_b(2,3,v_1,v_2)\geq s.
\)

We now extend the lower bound to arbitrary
real numbers $v_1,v_2\geq0$.
Choose positive rational sequences
\[
v_{1,j}\downarrow v_1,
\qquad
v_{2,j}\downarrow v_2.
\]
Since $v_{1,j}\geq v_1$ and $v_{2,j}\geq v_2$,
the inequalities defining
$R_b(2,3,v_{1,j},v_{2,j})$
are at least as restrictive as those defining
$R_b(2,3,v_1,v_2)$.
Then
\(R_b(2,3,v_{1,j},v_{2,j})
\subseteq R_b(2,3,v_1,v_2).
\)
The lower bound for positive rational
exponents gives
\[
\dimH R_b(2,3,v_1,v_2)
\geq D_{2,3}(v_{1,j},v_{2,j}).
\]
Letting $j\to\infty$ and using the continuity
of $D_{2,3}$, we obtain
\[
\dimH R_b(2,3,v_1,v_2)
\geq D_{2,3}(v_1,v_2).
\]
Together with the upper bound, this completes
the proof of Theorem~\ref{thm:exp}.

\section*{Acknowledgments}
The first author is partially supported by Beijing Natural Science Foundation (JQ26002).
The third author is partially supported by  Guangdong University
Young Innovative Talents Program Project (No. 2025KQNCX059).

\end{document}